\documentclass[a4paper,12pt]{scrartcl}

\usepackage[bookmarks,pagebackref]{hyperref}
\usepackage{amsrefs}
\usepackage[utf8]{inputenc}
\usepackage[english]{babel}
\usepackage{amssymb}
\usepackage{amsmath}
\usepackage{latexsym}
\usepackage{amsthm}
\usepackage{eucal}
\usepackage{bbm}
\usepackage{chngcntr}
\usepackage{apptools}
\usepackage{mathtools}
\usepackage{authblk}
\usepackage[pdflatex]{crop}
\usepackage{tikz}
\usepackage{tikz-cd}
\usepackage{microtype}
\usepackage{enumitem}
\usepackage{graphicx}
\usepackage{mathrsfs}
\usepackage[colorinlistoftodos]{todonotes}
\usepackage{yfonts}
\usepackage{tabularx}
\usepackage{color}

\allowdisplaybreaks

\setkomafont{disposition}{\normalfont\bfseries}
\newcommand{\auth}[0]{{Renan Assimos and J\"urgen Jost}}
\newcommand{\tit}[0]{{Harmonic maps from surfaces of arbitrary genus into spheres}}
\newcommand{\kw}[0]{{Bernstein theorem, minimal graph, harmonic map, Grassmannian, Gauss map, maximum principle}}

\hypersetup{
pdfauthor={\auth},%
pdftitle={\tit},%
colorlinks, linktocpage=true, pdfstartpage=1, pdfstartview=FitV,%
breaklinks=true, pdfpagemode=UseNone, pageanchor=true, pdfpagemode=UseOutlines,%
plainpages=false, bookmarksnumbered, bookmarksopen=true, bookmarksopenlevel=1,%
hypertexnames=true, pdfhighlight=/O,%
urlcolor=black, linkcolor=black, citecolor=black, 
}
\numberwithin{equation}{section}

\theoremstyle{plain}
\newtheorem{thm}{Theorem}[section]
\AtAppendix{\counterwithin{thm}{section}}
\newtheorem{defn}[thm]{Definition}
\newtheorem{lemma}[thm]{Lemma}

\newtheorem{notation}[thm]{Notation}

\newtheorem{rmk}[thm]{Remark}

\theoremstyle{definition}

\newtheorem{eg}[thm]{Example}

\newcommand{\Z}{\mathbb{Z}}

\newcommand{\C}{\mathbb{C}}
\newcommand{\R}{\mathbb{R}}

\DeclareMathOperator{\tr}{tr}

\let\originalleft\left
\let\originalright\right
\renewcommand{\left}{\mathopen{}\mathclose\bgroup\originalleft}
\renewcommand{\right}{\aftergroup\egroup\originalright}

\title{\tit}
\author{\auth\thanks{Correspondence: \href{mailto:assimos@mis.mpg.de}{assimos@mis.mpg.de}, \href{mailto:jjost@mis.mpg.de}{jjost@mis.mpg.de}}}
\affil{\small Max Planck Institute for Mathematics in the Sciences\\ Leipzig, Germany}
\date{}

\begin{document}

\maketitle

\begin{abstract}
We relate the existence problem of harmonic maps into $S^2$ to the convex geometry of $S^2$. On one hand, this allows us to construct new examples of harmonic maps of degree 0 from compact surfaces of arbitrary genus into $S^2$. On the other hand, we produce new example of regions that do not contain closed geodesics (that is, harmonic maps from $S^1$) but do contain images of harmonic maps from other domains. These regions can therefore not support a strictly convex function. Our construction builds upon an example of W. Kendall, and uses M. Struwe's heat flow approach for the existence of harmonic maps from surfaces.

\end{abstract}

\textbf{Keywords: }{Harmonic maps, harmonic map heat flow, maximum principle, convexity.}

\newpage

\tableofcontents

\section{Introduction}
M. Emery  \cite{Emery85} conjectured that a region in a Riemannian manifold that does not contain a closed geodesic supports a strictly convex function. W. Kendall  \cite{kendall92} gave a counterexample to that conjecture. Having such a counterexample is important to understand the relation between convexity of domains and convexity of functions in Riemannian geometry.

Here, we refine the analysis of that  counterexample and connect it with the existence theory of harmonic maps into $S^2$. The counterexample is $S^2$, equipped with its standard metric, minus three equally spaced subarcs of an equator of length $\pi/3$ each. That is, we subdivide the equator of $S^2$ into 6 arcs of equal length and remove every second of them. We construct nontrivial harmonic maps from compact Riemann surfaces whose image is compactly contained in that region. This excludes the existence of a strictly convex function on our region. The reason is, that the composition of a harmonic map with a convex function is a subharmonic function on the domain of the map, and therefore, since our domain is compact, has to be constant. But if both the map and the function are nontrivial, this is not possible. 

The main technical achievement of this paper is the construction of the harmonic maps. The existence theory of harmonic maps into surfaces was developed in important papers of J. Sacks and K. Uhlenbeck \cite{sacks81}, and L. Lemaire \cite{Lemaire78a}. Since then, various other existence schemes have been discovered (see for instance \cite{jost91}). While our construction will be more explicit than the general existence results, we also need arguments from the general theory. We use results from the heat flow approach of M. Struwe \cite{struwe85} \\

Let us now describe our methos in more technical terms. 
Let $(N,h)$ be a complete connected Riemannian  manifold and let $V\subsetneq N$ be an open connected subset of $N$. Suppose there exists a $C^2$-function $f:V\longrightarrow \R$ which is strictly geodesically convex. That is, 
\begin{equation*}
\left(f\circ c\right)^{\prime\prime}>0
\end{equation*}
for every geodesic $c:\left(-\epsilon,\epsilon\right)\longrightarrow V$. We say that such $V$ is a \textit{convex supporting domain}.

\begin{rmk}
While any subdomain of a convex supporting domain is convex supporting, we can not say much about their topology. But there are geometric constraints. For instance, for  $\epsilon>0$ there exists no strictly convex function supported in $\Omega_2^{\epsilon}:=\{x\in S^n\,|\,\, d(x,(0,0,1))\leq\frac{\pi}{2}+\epsilon\}$.   
\end{rmk}

\noindent Consider the following properties of an open connected subset $V\subsetneq N$.
\begin{equation*}
\begin{array}{l}
\textbf{(i)}\text{There exists a strictly convex function}\,\, f:V\longrightarrow \R\\

\textbf{(ii)}\,\,\text{There is no closed}\,\,M\,\,\text{and non-constant harmonic map}\,\,\phi:M\rightarrow N\,\,\text{such that}\,\,\phi(M)\subseteq V\\

\textbf{(iii)}\,\,\text{There is no}\,\,n\geq 1\,\,\text{and non-constant harmonic map}\,\,\phi:S^{n}\rightarrow N\,\,\text{such that}\,\,\phi(S^n)\subseteq V\\

\textbf{(iv)}\,\,\text{There is no  non-constant harmonic map}\,\,\phi:S^1\rightarrow N\,\,\text{such that}\,\,\phi(M)\subseteq V.
\end{array}
\end{equation*}

As already mentioned,  the maximum principle implies $\textbf{(i)}\implies \textbf{(ii)}$, and $\textbf{(ii)}\implies \textbf{(iii)}\implies \textbf{(iv)}$ are obvious. M. Emery \cite{Emery85}  conjectured that, if a subset $V\subseteq N$ admits no closed geodesics, that is, if it has property \textbf{(iv)}, then there exists a strictly convex function $f:V\rightarrow \R$; in other words,  \textbf{(iv)} should imply \textbf{(i)}. This conjecture was refuted by W. Kendall \cite{kendall92}, who gave an example of a subset of $S^2$ without closed geodesics. Kendall proves that one cannot define a strictly convex function in that subset of $S^2$. In the present paper we prove something slightly stronger regarding the existence of a non-constant harmonic map into the region constructed in \cite{kendall92}. Namely, we construct a non-constant harmonic into it. In other words, we prove that property \textbf{(iv)} does not imply \textbf{(ii)}. 

We use the harmonic map heat flow to construct a smooth harmonic map $u_{\infty}:(\Sigma_2,g)\longrightarrow S^2$, where $\Sigma_2$ is a closed Riemann surface of genus 2, and $u_{\infty}\subset \Omega\subset S^2$, where $\Omega:=S^2\setminus \left(\Gamma^1_{\epsilon},\Gamma^2_{\epsilon},\Gamma^3_{\epsilon}\right)$ is the region constructed in \cite{kendall92}, see Figure \ref{Fig5}. Our construction also works for Riemannian surfaces of higher genus, and it  gives a general method to construct harmonic maps from compact Riemann surfaces of any genus into symmetric surfaces. We first construct harmonic maps from graphs and then fatten the graphs to become Riemann surfaces. To construct the appropriate harmonic maps is the main technical achievement. In our construction, image symmetries under the harmonic map flow are preserved by the solution at time $t=\infty$.

What about other reverse implications for \textbf{(i)}-\textbf{(iv)}? In the last section of the present paper we show that the first reverse implication essentially holds following our work \cite{assimos2018}, and add remarks and strategies to decide whether the other implications might be true. 

A complete answer to all possible implications between \textbf{(i)}-\textbf{(iv)} would be a powerful result. In particular, since we have given abstract ways of showing that a region has property \textbf{(ii)} \cite{assimos2018}, this would yield a conceptual way of verifying whether a region does admit a strictly convex function. A similar question is raised by W. Kendall from a  martingale perspective \cite{kendall91}.

We thank Liu Lei and Wu Ruijun for helpful discussions. We thank Peter Topping for pointing out an imprecision in an earlier version of this paper. We thank Sharwin Rezagholi and Zachary Adams for helpful comments.

\section{Preliminaries}

\subsection{Harmonic maps} 

In this section, we introduce basic notions and properties of harmonic maps; references are \cite{xin2012}, \cite{lin2008} or \cite{Jost17}.

Let $(M,g)$ and $(N,h)$ be Riemannian manifolds without boundary, of dimension $m$ and $n$, respectively. By Nash's theorem there exists an isometric embedding $N\hookrightarrow \R^L$. 

\begin{defn}
	A map $u \in W^{1,2}(M,N)$ is called harmonic if and only if it is a critical point of the energy functional\\
	\begin{equation}\label{harm map def}
	E(u) := \frac{1}{2} \int_{M} \Vert du\Vert^{2}dvol_{g}
	\end{equation}
	where $\Vert .\Vert^{2} = \langle . ,.\rangle$ is the metric over the bundle $T^{*}M\otimes u^{-1}TN$ induced by $g$ and $h$.
\end{defn}

\noindent Recall that the Sobolev space $W^{1,2}(M,N)$ is defined as
\begin{align*}
W^{1,2}(M,N) =\hspace{0.2cm} &\Bigg\{v: M\longrightarrow \R^L; \hspace{0.2cm}\left|\left|v\right|\right|^{2}_{W^{1,2}(M)}=\int_{M}(\left|v\right|^2 + \Vert dv\Vert^{2})\hspace{0.1cm}dvol_g< +\infty \hspace{0.1cm}\text{and}\hspace{0.1cm}\\
&v(x)\in N \hspace{0.1cm}\text{for}\hspace{0.1cm}\text{almost}\hspace{0.1cm}\text{every}\hspace{0.1cm} x\in M \Bigg\}
\end{align*} 

\begin{rmk}
	If $M$ is a Riemann surface and $u \in W^{1,2}(M,N)$ is harmonic and $\phi: M\longrightarrow M$ is a conformal diffeomorphism, i.e. $\phi^* g=\lambda^2 g$ where $\lambda$ is a smooth function, we have, by \eqref{harm map def}, that
	\begin{equation}
	E(u\circ\phi)= \int_{M} \Vert du\Vert^{2}dvol_{g},
	\end{equation}
 that is, the energy of harmonic maps from surfaces is  conformally invariant. Therefore, on a surface, we need not specify a Riemannian metric, but only a conformal class, that is the structure of a Riemann surface, to define harmonic maps.
\end{rmk}

\noindent For $u\in W^{1,2}(M,N)$, define the map $du: \Gamma(TM) \longrightarrow u^{-1}(TN)$, given by $X \mapsto u_{*}X$. We denote by $\nabla du$ the gradient of $du$ over the induced bundle $T^*M\otimes u^{-1}(TN)$, that is, $\nabla du$ satisfies $\nabla_Y du\in\Gamma (T^*M\otimes u^{-1}(TN))$, for each $Y\in\Gamma(TM)$.

\begin{defn}[Second fundamental form]
	The second fundamental form of the map $u:M\longrightarrow N$ is the map defined by
	\begin{equation}\label{sec fund form 1}
	B_{X Y}(u)=(\nabla_X du)(Y)\in \Gamma(u^1 TN).
	\end{equation}
	$B$ is bilinear and symmetric in $X,Y\in \Gamma(TM)$. It can also be seen  as 
	\begin{equation*}
	B(u)\in\Gamma(Hom(TM\odot TM, u^{-1} TN)).
	\end{equation*}
\end{defn}

\begin{defn}
	The tension field of a map $u:M\longrightarrow N$ is the trace of the second fundamental form
	\begin{equation}
	\tau(u)=B_{e_i e_i}(u)=(\nabla_{e_i} du)(e_i)
	\end{equation}
	seen as a cross-section of the bundle $u^{-1} TN$.
\end{defn}

\noindent Taking a $C^1$ variation $\big(u(.,t)\big)_{\|t\|<\epsilon}$ yields
\begin{equation}\label{first var of energy}
	\frac{d}{dt} E(u(.,t)) = -\left\langle\tau(u(x,t)),\frac{\partial u}{\partial t}\right\rangle_{L^2}.
\end{equation}

\noindent By \eqref{first var of energy}, the Euler-Lagrange equations for the energy functional are
\begin{equation}\label{harmonic map equation}
\tau(u)= \Bigg(\Delta_{g} u^{\alpha} + g^{ij}\Gamma_{\beta\gamma}^{\alpha}\frac{\partial u^{\beta}}{\partial x^i}\frac{\partial u^{\gamma}}{\partial x^j}h_{\beta\gamma}\Bigg)\frac{\partial}{\partial u^{\alpha}} = 0,
\end{equation}
where the $\Gamma_{\beta\gamma}^{\alpha}$ denote the Christoffel symbols of $N$.

\noindent Hence we can equivalently define (weakly) harmonic maps as maps that (weakly) satisfy the \textit{harmonic map equation} \eqref{harmonic map equation}.

\begin{defn}[Totally geodesic maps]
	A map $u\in W^{1,2}(M,N)$ is called totally geodesic if and only if its second fundamental form identically vanishes, i.e. $B(u)\equiv 0$. 
\end{defn}

\begin{lemma}[Composition formulas]
	Let $u:M\longrightarrow N$ and $f:N\longrightarrow P$, where $(P,i)$ is another Riemannian manifold. Then 
	\begin{align}\label{composition formulas}
	\nabla d(f\circ u) = df\circ\nabla du + \nabla df(du,du),\\
	\tau(f\circ u) = df\circ\tau(u) + \tr\nabla df(du,du). 
	\end{align} 
\end{lemma}

\begin{eg}[L. Lemaire \cite{lemaire78}]
	
	For every integers $p$ and $\mathcal{D}$ such that $|\mathcal{D}\leq p-1|$, there exists a Riemann surface of genus $p$ and a harmonic nonholomorphic map of degree $\mathcal{D}$ from that surface to the sphere.
	
	\noindent For the case of $\mathcal{D}=0$, one considers a 1-form as a harmonic map  $\tilde{\phi}:\Sigma_g\longrightarrow S^1$ and a totally geodesic embedding $S^1\hookrightarrow S^2$. The composition $\phi := \,(\text{tot.geod})\, \circ\, \tilde{\phi} $ is harmonic. Here, the image $\phi(\Sigma_g)\subset S^2$ is contained in the image of a closed geodesic of $S^2$, and therefore, is geometrically trivial. The main construction of this paper  yields harmonic maps from  Riemann surfaces of even genus into $S^2$ that do not contain closed geodesics in their images.
\end{eg}

\subsection{The harmonic map flow}
Let $M$ be a compact Riemann surface and $N$ a compact Riemannian manifold. 
M. Struwe \cite{struwe85} (see also Chang \cite{Chang89}) showed that the heat flow
\begin{equation}\label{hmf eq}
\begin{cases}
\partial_{t}u(x,t) = \tau\big(u(x,t)\big)\\
u(\hspace{3pt}\cdot\hspace{3pt},0) = u_0.
\end{cases}
\end{equation}
has a global weak solution. 
By \eqref{first var of energy}, we have
\begin{eqnarray}\label{energy decreases along the flow}
\frac{d}{dt}E(u(x,t)) &=& -\int_{M}\big\langle\tau\big(u(x,t)\big),\partial_{t}u(x,t)\big\rangle_{u^{\ast}TN}dvol_{g}\nonumber\\
&=& -\int_{M}\big|\tau(u(x,t))\big|^2dvol_{g},
\end{eqnarray}
and since $\big|\tau(u(x,t))\big|^2 \geq 0$, the energy does not increase along the flow. 

M. Struwe's solution  is smooth with the possible exception of finitely many singular points.  A. Freire \cite{freire95uniqueness1} proved uniqueness among the maps for which $E(\,\cdot\,,t)$ is decreasing in $t$, see P. Topping \cite{topping2002,topping2010}. More precisely, we have the following.

\begin{thm}[M. Struwe, A. Freire]\label{Struwe 1}
	For any initial value $u_0\in W^{1,2}(M,N)$, there exists a weak solution $u$ of the equation~\eqref{hmf eq} in $W^{1,2}(M\times[0,+\infty),N)$, and this solution is unique among the maps for which $E(\,\cdot\,)$ is decreasing in $t$. Moreover, in $M\times[0,+\infty)$, $u$ is smooth with the exception of finitely many points.
	
\end{thm}

\noindent At each such singular point, we encounter  the bubbling phenomenon.


\noindent   Struwe's results on the harmonic map flow will be  the main tools to construct the desired harmonic map.

\section{The main example}

\subsection{The construction of the Riemann surface}

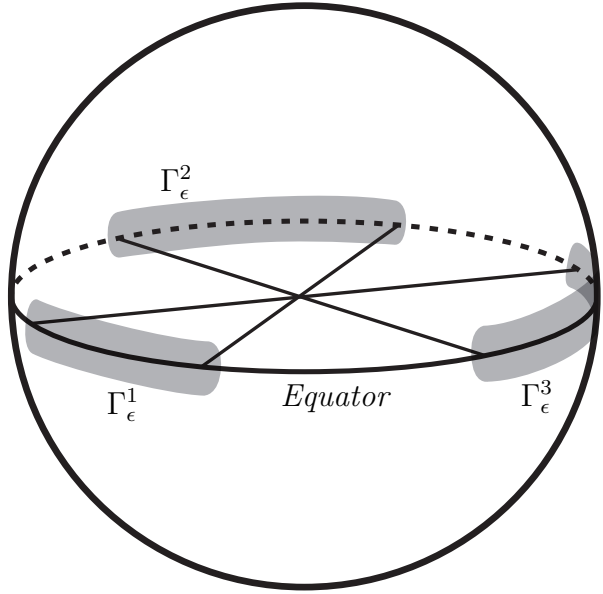
\begin{figure}
	\begin{picture}(200,200)
	\put(85,0){\includegraphics[width=0.6\linewidth]{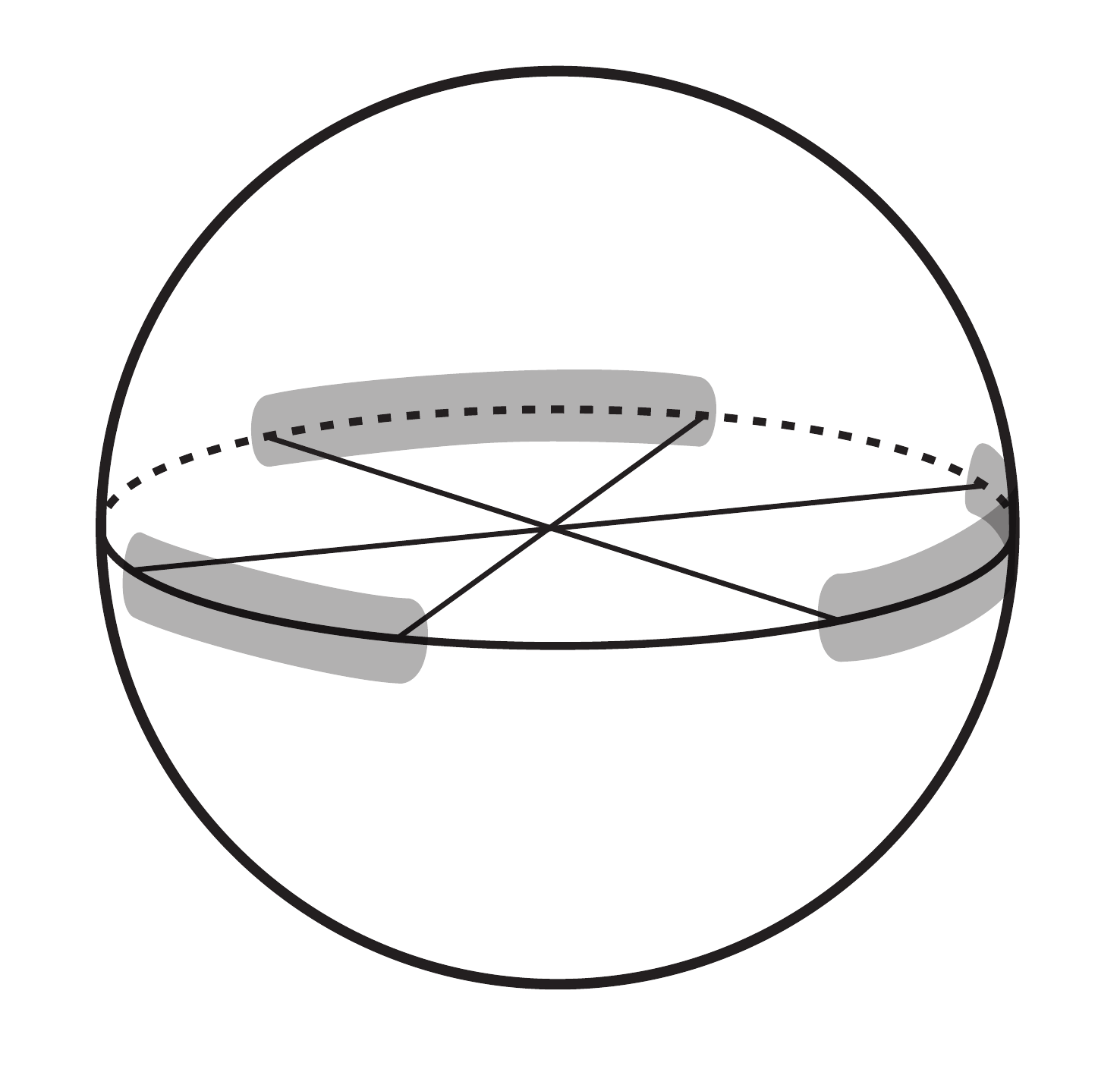}}
	\put(145,85){$\Gamma^1_{\epsilon}$}
	\put(165,169){$\Gamma^2_{\epsilon}$}
	\put(300,88){$\Gamma^3_{\epsilon}$}
	\put(210,88){$\textit{Equator}$}
	\end{picture}
	\caption{A region in $S^2$ without closed geodesics.}
	\label{Fig6}
\end{figure}

We start by considering a region on $S^{2}$ described by Figure~\ref{Fig6}.

\begin{notation}
	The `$\text{Equator}$' is the embedding $i:S^{1}\hookrightarrow S^{2}$ with $i(S^1)=\{(x_{1},x_{2},x_{3})\in S^{2}; \hspace{3pt}x_{3}=0\}$.
\end{notation}

\noindent Dividing the $\textit{Equator}$ into 6 equal pieces, we can pick $\Gamma^1$, $\Gamma^2$, $\Gamma^3$ as three non-intersecting pieces, as in Figure~\ref{Fig6}. For a given $\epsilon>0$, consider the sets $\Gamma^i_{\epsilon}:=\{x\in S^2|\,\, d_{\mathring{g}}(x,\Gamma^i)<\epsilon\}$, for each $i=1,2,3$. Define $\Omega_{\epsilon}:=S^2\setminus\left(\Gamma^1_{\epsilon}\cup\Gamma^2_{\epsilon}\cup\Gamma^3_{\epsilon}\right)$.

For a point $p\in  \textit{Equator}\setminus(\Gamma^{1}\cup\Gamma^{2}\cup\Gamma^{3})$, it is obvious that its antipodal point $A(p)$ fulfills $A(p)\in (\Gamma^{1}\cup\Gamma^{2}\cup\Gamma^{3})$. Therefore $\Omega_{\epsilon}$ does not admit closed geodesics. This is analogous to W. Kendall's propeller \cite{kendall92}.

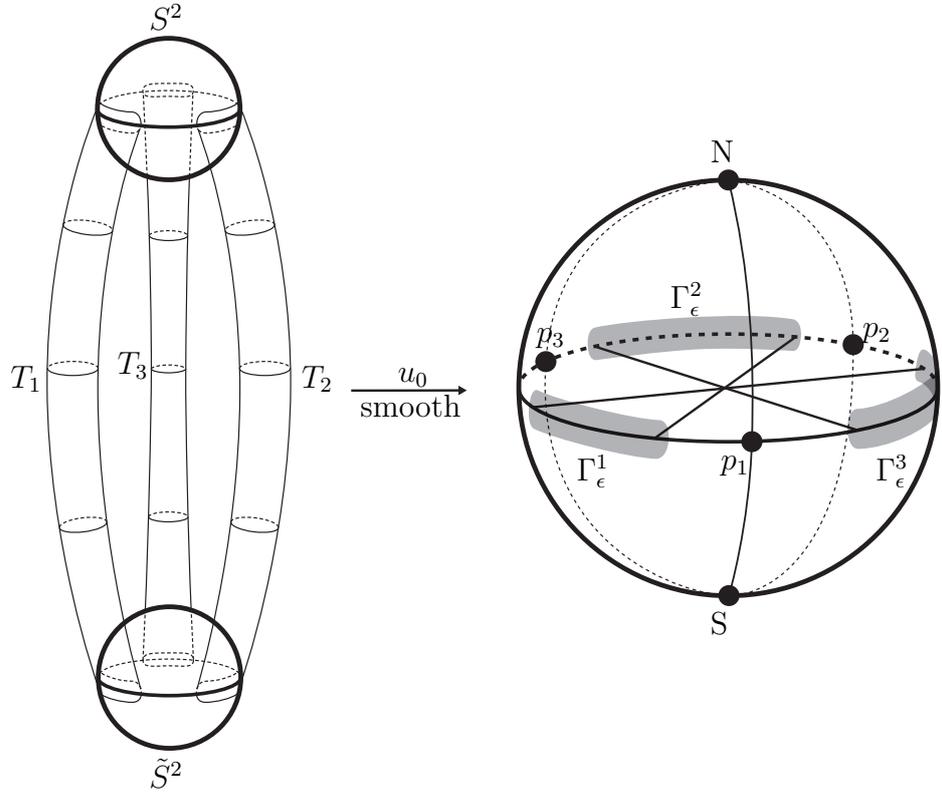
\begin{figure}
	\begin{picture}(100,260)
	\put(70,0){\includegraphics[width=0.85\linewidth]{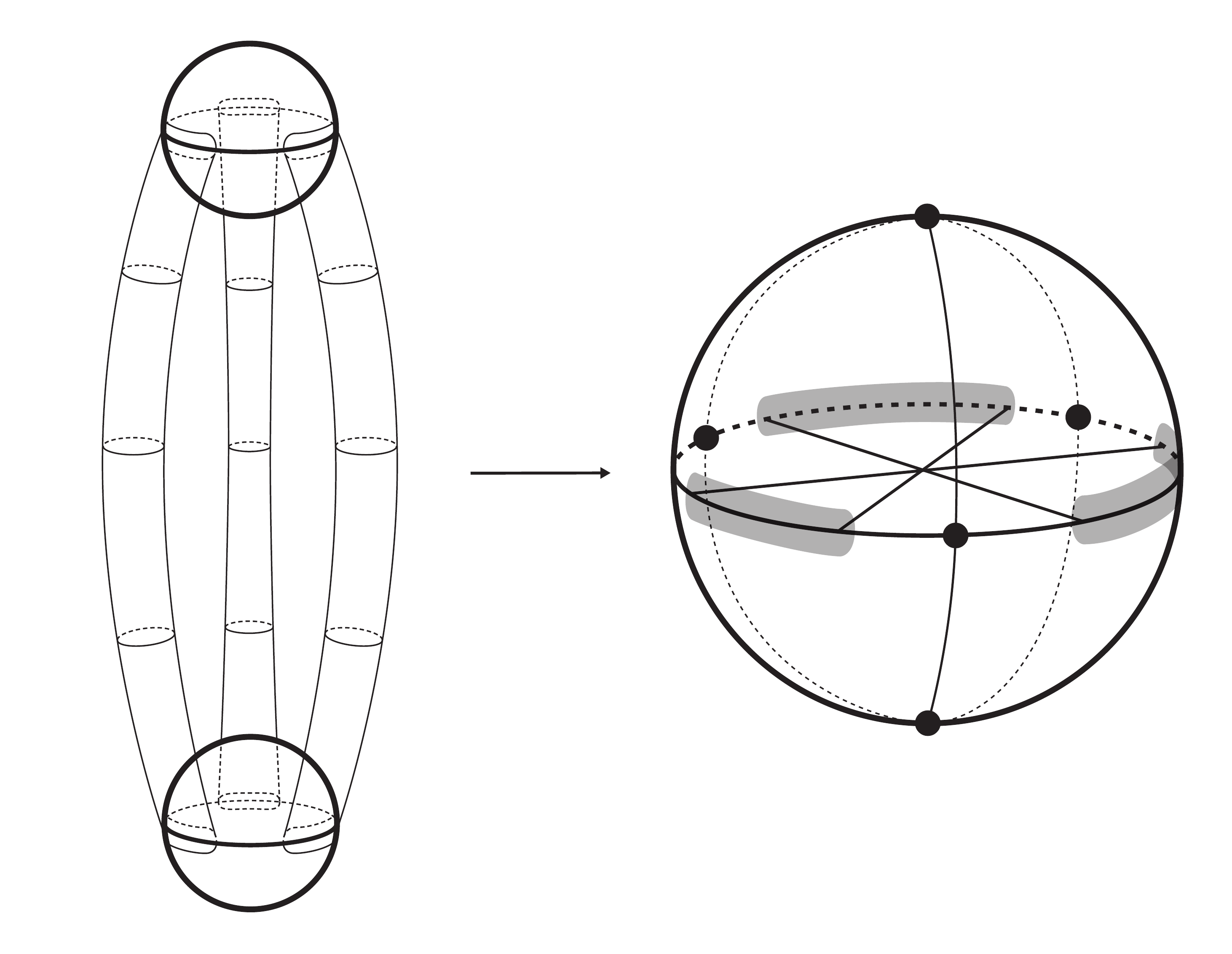}}
	\put(233,152){$u_0$}
	\put(219,139){smooth}
	\put(128,155){$T_3$}
	\put(88,150){$T_1$}
	\put(197,150){$T_2$}
	\put(140,285){$S^2$}
	\put(140,-1){$\tilde{S}^2$}
	\put(350,235){N}
	\put(350,57){S}
	\put(300,115){$\Gamma^1_{\epsilon}$}
	\put(412,115){$\Gamma^3_{\epsilon}$}
	\put(335,180){$\Gamma^2_{\epsilon}$}
	\put(285,167){$p_3$}
	\put(354,119){$p_1$}
	\put(407,169){$p_2$}
	\end{picture}
	\caption{The Riemann surface $\Sigma_2$ and the map $u_0$}
	\label{Fig5}
\end{figure}

We define a compact Riemann surface $\Sigma_{2}$ as follows, see Figure~\ref{Fig5}.

Consider the poles $S = (0,0,-1)$ and $N=(0,0,1)$. Let $\tilde{S}^{2}$ denote another sphere such that $\tilde{N}$, $\tilde{S}$ belong to the axis $x_{3}$ of $R^3$ and $\tilde{S}^{2}$ is a reflection of $S^{2}$ centered at a point of the axis. Denote by $\text{d} > 0$ the distance between $\tilde{N}$ and $S$. In $\tilde{S}^2$, consider also the segments $\tilde{\Gamma}^{i}$ and the sets $\tilde{\Gamma}^i_{\epsilon}$ as in $S^2$.

We connect the two spheres by three tubes (cylinders), called $T^1$, $T^2$, $T^3$ such that $\partial T^i = \Gamma^i_{\epsilon}\cup\tilde{\Gamma}^i_{\epsilon}$, i.e., the tubes connect the sets $\Gamma^i_{\epsilon}$ and $\tilde{\Gamma}^i_{\epsilon}$ pairwise and smoothly.

\noindent 
The following two symmetries will be important.
\begin{enumerate}
	\item A $\Z_{3}$ symmetry around the $x_{3}$ axis is preserved, i.e. the tubes and the circles are all equal.
	\item A $\Z_{2}$ symmetry around the midpoint $\tilde{N}$ and $S$ is preserved.
\end{enumerate} 

\subsection{The initial condition $u_0$}

The procedure on the previous subsection generates a genus 2 compact Riemann surface $\Sigma_{2}$, and we consider the induced metric from $\R^3$ as its Riemannian metric. The next step is to define the smooth initial condition between $\Sigma_{2}$ and $S^2$, and explore the properties of the harmonic map flow.

Let $u_0\in C^{\infty}\big(\Sigma_2,S^2\big)$ be the map given by 
\begin{eqnarray}
&u_0: x  & \mapsto \begin{cases}
\begin{array}{ccl}
N& \hspace{3pt}\text{if}&\hspace{3pt} x\in S^2\setminus(C^1\cup C^2\cup C^3) \\
S& \hspace{3pt}\text{if}&\hspace{3pt} x\in \tilde{S}^2\setminus(C^1\cup C^2\cup C^3) \\
(*)& \hspace{3pt}\text{if}&\hspace{3pt} x\in (T^1\cup T^2\cup T^3).
\end{array}
\end{cases}
\end{eqnarray}

\noindent Where $(*)$ is defined as follows. Consider the geodesic $\gamma$ connecting the south and the north poles in $S^2$ given below.
\begin{equation*}
\begin{array}{rcl}
\gamma:&[-R,+R]& \longrightarrow S^2\\
& t& \longmapsto \Bigg(\sin\left(\frac{\pi(t+R)}{2R}\right), 0 , -\cos\left(\frac{\pi(t+R)}{2R}\right)\Bigg),
\end{array}
\end{equation*}
where $2R > 0$ is the height of the tubes $T^i$ in our Riemann surface $\Sigma_2$. Now, given $\theta\in\left[0,2\pi\right)$, 
\begin{equation*}
\begin{array}{rcl}
\Gamma_{\theta}: [-R,R]& \longrightarrow& T^i\\
t &   \longmapsto & (r\cos\theta, r\sin\theta, t),
\end{array}
\end{equation*}
where $r$ is the radius of $T^i$. 

Up to the diffeomorphism that makes the cylinder straight, define $u_0|_{T^i}$ as a map in $C^{\infty}(\Sigma_2,S^2)$ such that
\begin{equation}\label{def of u0}
u_0\circ\Gamma_{\theta}(t)=\gamma(t)\,\,\text{for}\,\text{every}\,\,\theta\in\left[0,2\pi\right).
\end{equation}

The image of $\Sigma_2$ under $u_0$ preserves the $\Z_{2}$ and $\Z_{3}$ symmetries of $\Sigma_2$, and the harmonic map flow will preserve these symmetries of the initial condition.

\begin{lemma}{(Symmetry property of the flow)}
	Let $\varphi$ be the $\Z_{3}$ action on $\Sigma_{2}$ given by the rotation by $\frac{3\pi}{2}$ over the axis connecting the poles. Let $\psi$ be the `same' $\Z_3$ action over the sphere $S^2$. Obviously $u_0\circ\varphi = \psi\circ u_0$, hence we have
	
	\begin{equation}
	\begin{cases}
	\partial_{t}\big(\psi\circ u\big) = \tau\big(\psi\circ u\big)\\
	\psi\circ u(\hspace{3pt}\cdot\hspace{3pt},t) = u_0\circ\varphi
	\end{cases}
	\end{equation}
	or equivalently
	\begin{equation}
	\begin{cases}
	\partial_{t}\big(u\circ\varphi\big) = \tau\big(u\circ\varphi\big)\\
	u\circ\varphi (\hspace{3pt}\cdot\hspace{3pt},t) = \psi\circ u_0.
	\end{cases}
	\end{equation}
	
	\noindent Moreover, if $u$ is a the unique Struwe solution for \eqref{hmf eq} with initial condition \eqref{def of u0}, as $\psi$ is totally geodesic, we have
	\begin{eqnarray*}
		\partial_t(\psi\circ u) - \tau(\psi\circ u) &=&  \text{d}\psi(u(\hspace{3pt}\cdot\hspace{3pt},t))\partial_t u - \text{d}\psi(u(\hspace{3pt}\cdot\hspace{3pt},t))\tau(u)\\
		&=& \text{d}\psi(u(\hspace{3pt}\cdot\hspace{3pt},t))\big(\partial_t u - \tau(u)\big)=0,
	\end{eqnarray*}
	
	\noindent therefore $\psi\circ u$ is a solution to the same problem.
\end{lemma}

 The same argument applies to the $\Z_2$ action $(a,b,c)\mapsto (a,b,-c)$. Therefore the harmonic map $u_{\infty}:\Sigma_2\rightarrow S^2$ given by Struwe's solution of \eqref{hmf eq} with initial condition $u_0$ is $\Z_3$ and $\Z_2$ equivariant with respect to the actions defined above.

 Another important remark about the initial condition $u_0$ is that we can control its energy by changing the length and radius of the tubes $T^i$.

\begin{lemma}[A control on the energy of $u_0$]\label{control on the energy}
	
	Since $\|\dot{\gamma}(t)\|=\frac{\pi}{2R}$, we have
	\begin{equation*}
	0 < \| \text{d}u_0 (p)\|\leq\frac{\pi}{2R}\hspace{7pt}
	\end{equation*}
	and therefore
	\begin{eqnarray*}
		0 < E(u_0) \leq& \displaystyle \frac{1}{2}\int_{T^i}\left( \frac{\pi}{2R}\right)^2dvol_{T^i}\\ 
		=& \displaystyle \frac{\pi^2}{8R^2}\int_{T^i}dvol_{T^i}\\
		=& \displaystyle \frac{\pi^3 r^2}{4R}.
	\end{eqnarray*}

	\noindent By making the tubes connecting $S^2$ and $\tilde{S^2}$ in $\Sigma_2$ thinner and longer, we can make the energy of $u_0$ arbitrarily small. More precisely, given any $\epsilon>0$, we can pick $r\in(0,1)$ and $R\in(1,+\infty)$ such that $E(u_0)\leq\epsilon$. By Equation~\eqref{energy decreases along the flow}, we know that the energy decreases along the flow and therefore we have a control on the energy of $u_t$ for every $t\in \R_+$.
\end{lemma}

\subsection{Controlling the image of $u_{\infty}$}

Theorem~\ref{Struwe 1} above roughly tells us that, given any initial condition $u_0\in W^{1,2}(M,N)$, there exist a smooth solution to the harmonic map flow with the exception of some finite points on which the energy is controlled. In other words, we have a harmonic map
\begin{equation*}
u_{\infty}: M\backslash\{q_{1},...,q_{l} \} \rightarrow N
\end{equation*}
and around the singularities this map can be extended to a harmonic sphere $h: S^2 \rightarrow N$.

Since each harmonic two-sphere which bubbles out carries at least the energy $4\pi$ (the area of the target sphere $S^2$)  and every energy loss during the flow is due to the formation of a bubble (see \cite{qingtian97}),  we avoid the formation of bubbles by taking $\Sigma_2$ as the compact Riemann surface with tubes of length $R\in(1,+\infty)$ and radius $r\in(0,1)$, such that for a given $\epsilon_0>0$ we have $E(u_0)\leq\epsilon_0 << 4\pi$.

With this initial condition $u_0$, we have a unique global smooth solution $u$ of the harmonic map flow \eqref{hmf eq}. This gives us the smooth harmonic map
\begin{equation}
u_{\infty}:\Sigma_2 \rightarrow S^2.
\end{equation}

It remains to  prove that $u_{\infty}(\Sigma_2)\subset S^2\setminus (\Gamma^{1}\cup\Gamma^{2}\cup\Gamma^{3})$.

\noindent To do so, we need the following lemma due to Courant, see for instance \cite{jost91} for a proof. 

\begin{lemma}[Courant-Lebesgue]\label{c-l}
	Let $f\in W^{1,2}(D,\R^d)$, $E(f)\leq C$, $\delta<1$ and $p\in D=\{(x,y)\in\C;\, x^2+y^2=1\}$. Then there exists some $r\in(\delta,\sqrt{\delta})$ for which $\displaystyle f|_{\partial B(p,r)\cap D}$ is absolutely continuous and
	\begin{equation}
	\displaystyle\left|f(x_1)-f(x_2)\right|\leq (8\pi C)^{\frac{1}{2}}\left(\log \frac{1}{\delta}\right)^{-\frac{1}{2}}
	\end{equation}
	\noindent for all $x_1,x_2\in\partial B(p,r)\cap D$.
\end{lemma}

The $\Z_{2}$ symmetry of $u_{\infty}(\Sigma_2)$ implies that there are two antipodal points, called $N$ and $S$, on the image of $u_{\infty}(\Sigma_2)$. By the $\Z^3$ symmetry and the fact that the image is connected, there exist three points $p_i\in T_{i}$, $i=1,2,3$, each of them $\Z_2$-invariant, with $u_{\infty}(p_i)\in \textit{Equator}$. Moreover, $\varphi(u_{\infty}(p_1))=u_{\infty}(p_2)$, $\varphi(u_{\infty}(p_2))=u_{\infty}(p_3)$ and $\varphi(u_{\infty}(p_3))=u_{\infty}(p_1)$. With the help of  the Courant-Lebesgue lemma \ref{c-l}, we shall now see that that a small neighborhood of the tube around $p_i$ is mapped into a small neighborhood of $u_{\infty}(p_i)$ of controlled size. In fact, since $E(u_\infty)\leq E(u_0)<\frac{r^2}{R}$, where $r$ is the radius of the tubes and $R$ their heights, we start taking $r\in(0,1)$ and $R\in(1,+\infty)$ such that $u_{\infty}$ has no bubbles and $r\pi<<\pi/6$. 

\noindent Therefore, taking $\delta<r^2$ the Courant-Lebesgue lemma implies that there exists some $s\in(r^2,r)$ such that
\begin{equation}\label{courant lebesgue distance}
|u_{\infty}(x_1)-u_{\infty}(x_2)|\leq \Bigg(\frac{8\pi.r^2}{R}\Bigg)^{\frac{1}{2}}\Bigg(\log\Big(\frac{1}{\delta}\Big)\Bigg)^{-\frac{1}{2}}
\end{equation}
for every $x_1,x_2\in\partial B_{g}(p_i,s)\cap B_g(p_i,1)$. But, since $\delta<r^2$, we get
\begin{equation*}
\begin{array}{rcl}
\frac{1}{\delta}>\frac{1}{s}& \Rightarrow &\log\left(\frac{1}{\delta}\right) > \log\left(\frac{1}{s}\right)\\ 
& &\\
& \Rightarrow & \left(\log\left(\frac{1}{s}\right)\right)^{-\frac{1}{2}} > \left(\log\left(\frac{1}{\delta}\right)\right)^{-\frac{1}{2}}.\\
\end{array}
\end{equation*}
Therefore, by Equation~\eqref{courant lebesgue distance},
\begin{equation*}
\begin{array}{rcl}
|u_{\infty}(x_1)-u_{\infty}(x_2)|&\leq & \left(\frac{8\pi.s}{R}\right)^{\frac{1}{2}}\left(\log\left(\frac{1}{s}\right)\right)^{-\frac{1}{2}}\\
& & \\
& = & \frac{\sqrt{s}}{\left(\log\left(\frac{1}{s}\right)\right)^{\frac{1}{2}}}.\sqrt{\frac{8\pi}{R}}<<\frac{\pi}{6}.
\end{array}
\end{equation*}

In particular, if we call $S^1_i$ the circle given by the intersection of the $x,y$-plane in $\R^3$ with $T_i$, the above argument shows that the image of a $\delta$-neighborhood of $S^1_i$ under $u_{\infty}$ is contained in a neighborhood of $u_{\infty}(p_i)$ that does not intersect any of the removed sets $\Gamma^i_{\epsilon}$, because, since $p_i$ is invariant under the $\Z_2$ action, so is $S^1_i$.

We claim that for any point $a\in T_i$ such that $a\notin B(S^1_i,\delta):=\{x\in T_i\,|\,\, d(x,S^1_i)\leq\delta\}$, it follows that $u_{\infty}(a)\notin \textit{Equator}$. 

Obviously $u_{\infty}(\Sigma_2)\subsetneq S^2$, since $Area(u_{\infty})\leq E(u_{\infty})<\frac{r^2}{R}<<1$. Since $\Sigma_2$ is compact, $u_{\infty}(\Sigma_2)$ is compact.

If $\partial u_{\infty}(\Sigma_2)=\emptyset$, then $u_{\infty}(T_i)$ is totally geodesic and we have that $u_{\infty}(a)\notin \textit{Equator}$ if $a\in T_i\setminus B(S^1_i,\delta)$. If $\partial u_{\infty}(\Sigma_2)\neq\emptyset$ we argue by contradiction. Suppose there exists $a\in T_{i}\setminus B(S^1_i,\delta)$ such that $u_{\infty}(a)\in\partial u_{\infty}(\Sigma_2)$ and $u_{\infty}(a)\in \textit{Equator}$. For $|t|<\epsilon$, consider the $C^1$-variation
\begin{equation}
\left(u_{\infty}\right)_t(x) = u_{\infty}(x) + t\eta(x),
\end{equation}
where $\eta:\Sigma_2\longrightarrow S^2$ is smooth, $\eta\equiv 0$ outside a neighborhood of $a$ and $\left(u_{\infty}\right)_t(a)\in \textit{Hemisphere}_N$ (where we are assuming without loss of generality that $d(a,S^2)<d(a,\tilde{S}^2)$). This variation $\left(u_{\infty}\right)_t(x)$ clearly satisfies $E\left(\left(u_{\infty}\right)_t\right)<E\left(u_{\infty}\right)$, but this contradicts the fact that $u_{\infty}$ is a smooth harmonic map. Therefore such a point $a\in T_i\setminus B(S^1_i,\delta)$ cannot be mapped on the $\textit{Equator}$.

\noindent This concludes the proof that $u_{\infty}(\Sigma_2)\subset S^2\setminus(\Gamma^1_{\epsilon}\cup\Gamma^2_{\epsilon}\cup\Gamma^3_{\epsilon})$.

To summarize, we have shown that the solution $u_{\infty}$ of the harmonic map flow

\begin{equation*}
\begin{cases}
\partial_{t}u(x,t) = \tau\big(u(x,t)\big)\\
u(\hspace{3pt}.\hspace{3pt},t = 0) = u_0
\end{cases}
\end{equation*}
where $u_0:\Sigma_2\rightarrow S^2$ is given by Definition~\ref{def of u0} is a harmonic map from the compact genus two Riemann surface $\Sigma_2$ into the sphere $S^2$ and the image $u_{\infty}(\Sigma_2)$ does not contain closed geodesics.

\begin{rmk}
	By the same argument, we could construct harmonic maps from compact Riemann surfaces of genus $2p$ for any $p>1$, replacing the $\Z_3$ symmetry  by a $\Z_{2p+1}$ symmetry. Moreover, if we replace the target by another $2$-dimensional surface with the appropriate symmetries, we could try to construct a $u_0$-type initial condition respecting the necessary symmetries.
\end{rmk}

\section{Remarks and open questions}

As  pointed out in the introduction, it is an interesting question whether the implications between properties $(i)$ to $(iv)$ of a given subset $V\subset (N,h)$ hold true. Here, we make a couple of remarks on this problem to illustrate some applications and possible consequences.

We start by considering the equivalence between properties  $(i)$ and $(ii)$. Since $(i)$ implies $(ii)$ by the maximum principle, the question is whether $(ii)$ implies $(i)$. More precisely, we are assuming that a subset $V\subseteq (N,h)$ admits no image of non-constant harmonic maps $\phi:(M,g)\longrightarrow (N,h)$, where $M$ is a closed Riemannian manifold, and we ask if this implies the existence of a strictly convex function $f:V\longrightarrow \R$.

We start by remarking that, if one takes a geodesic ball $B(p,r)$ in a complete manifold $N$ such that $r$ is smaller than the convexity radius of $N$ at $p$, then $\partial B(p,r)$ is a hypersurface of $N$ with definite second fundamental form for every point $q\in\partial B(p,r)$. 	

A partial answer to the above question was obtained in \cite{assimos2018} where the following theorem was proven.

\begin{thm}\label{corollary of SMP}
	Let $(N,h)$ be a complete Riemannian manifold and $\gamma:[a,b]\longrightarrow N$ a smooth embedded curve. Consider a smooth function $r:[a,b]\longrightarrow \R_+$ and a region
	\begin{equation}\mathcal{R}:=\bigcup_{t\in[a,b]}B(\Gamma(t),r(t)),\end{equation}
	where $B(\cdot,\cdot)$ is the geodesic ball and $r(t)$ is smaller than the convexity radius of N for any $t$.
	If, for each $t_{0}\in (a,b)$, the set $\mathcal{R} \backslash B(\gamma(t_{0}),r(t_{0}))$ is the union of two disjoint connected sets, namely the connected component of $\gamma(a)$ and the one of $\gamma(b)$, then there exists no compact manifold $(M,g)$ and non-constant harmonic map $\phi:M\longrightarrow N$ such that $\phi(M)\subseteq\mathcal{R}$. 
\end{thm}

According to the above definitions,  the region $\mathcal{R}$ of theorem~\ref{corollary of SMP} has property $(ii)$. As a direct corollary of the proof, which is based on an application of a maximum principle, we conclude that, if a subset $V$ of $N$ admits a sweep-out by convex hypersurfaces with the additional property that the absence of each leave of this sweep-out divides $V$ in two connected components like the region $\mathcal{R}$ above, then $V$ has property $(ii)$.

This corollary gives us a geometric way of checking whether  a certain subset of a manifold has property $(ii)$. If  $f:V\longrightarrow \R$ is strictly convex, then the level sets of $f$ naturally give us the desired sweepout in \ref{corollary of SMP}. On the other hand, a function $f:V\longrightarrow \R$ whose level sets yield a convex sweepout, is not necessarily strictly convex. We may have to reparametrize $f$ to make it  strictly convex: It is enough to use the parameter $t\in[a,b]$ of the 1-parameter family of convex leaves in the above theorem and control the growth of the gradient of $f$ while it walks through the convex leaves of the sweep-out.

The above argument is a strategy for the proof that when $V$ satisfies the properties of the region $\mathcal{R}$ in theorem~\ref{corollary of SMP}, then $V$ has property $(i)$. In other words, this gives us a method to obtain strictly convex functions on some subset $V$ of $(N,h)$ based on its geometry.

Another question is whether $(iv)$ implies $(iii)$. Suppose $V$ does not have closed geodesics. Does that imply that there are no harmonic maps $\phi:(S^k,\mathring{g})\longrightarrow (N,h)$ with $\phi(S^k)\subset V$? In particular, does the absence of closed geodesics imply no bubbles in $V$?

In the previous sections we have proven that property $(iv)$ does not imply property $(ii)$. It may be that with similar techniques, and a  symmetric space different from $S^2$ as target, one can build a counterexample for the implication  $(iv) \implies (iii)$ as well, allowing the energy of the initial condition to be big enough to form a bubble, while still controlling the image of the final map preserving symmetries, as in section 2.       

\bibliographystyle{alpha}
\bibliography{bernstein-harmonic}
\addcontentsline{toc}{section}{\bibname}

\end{document}